\documentclass{amsart}
\usepackage{amssymb,amsmath}
\usepackage{graphicx}
\usepackage[enableskew]{youngtab}
\usepackage{young}
\usepackage{tikz,pgf}

\newtheorem{theorem}{Theorem}

\newtheorem{conjecture}{Conjecture}

\begin{document}

\title{A bijection for partitions simultaneously $s$-regular and $t$-distinct}
\author{William J. Keith}
\maketitle

\begin{abstract}
In this note a bijection is constructed between the set of partitions of $n$ simultaneously $s$-regular and $t$-distinct, and those simultaneously $t$-regular and $s$-distinct.  Some implications of the map are discussed.  As a generalized version of Glaisher's bijection, the map may be widely useful in other partition combinatorics.  A previous conjecture concerning iterations of Glaisher's bijection is given a counterexample.
\end{abstract}

\section{Introduction}

A partition of $n$ is a nonincreasing sequence of positive integers $(\lambda_1, \lambda_2, \dots , \lambda_k)$ which sums to $n$.  We write $\lambda \vdash n$ or $\vert \lambda \vert = n$ if $\lambda_1 + \dots + \lambda_k = n$. Their study was initiated by Euler in \cite{Euler}, who proved the usual first theorem learned by a student of the area, namely

\begin{theorem} The number of partitions of $n$ in which all parts are odd equals the number of partitions of $n$ in which parts are distinct.
\end{theorem}

Euler's proof was by equating the two generating functions: $$\prod_{k=1}^\infty (1+q^k) = \prod_{k=1}^\infty \frac{1-q^{2k}}{1-q^{k}} = \prod_{k=1}^\infty \frac{1}{1-q^{2k-1}}.$$  The leftmost expression is most clearly the generating function for partitions where no part size repeats, and the rightmost expression is more obviously the generating function for partitions in which only odd parts are allowed.

The theorem was proved by a hands-on combinatorial mapping found by J. J. Sylvester, and then generalized to all moduli by a more general mapping given by his student Glaisher:

\begin{theorem} The number of partitions of $n$ in which no part is divisible by $m$ equals the number of partitions of $n$ in which parts appear fewer than $m$ times.
\end{theorem}

The generating functions for both are equal to $$\prod_{k=1}^\infty \frac{1-q^{mk}}{1-q^{k}}.$$

Henceforth we shall call the first type of partition $m$-\emph{regular}, and the second type $m$-\emph{distinct}.  (It should be noted that the term $m$-distinct also appears in the literature denoting partitions with parts differing by at least $m$.)

Glaisher's map is now a fundamental tool in the combinatorial theory of partitions.  It is as follows: if in an $m$-distinct partition part $j m^k$ appears $a_{j,k}$ times, $m \nmid j$, $a_{j,k} < m$, then write $a_{j,k} m^k$ appearances of $j$.  The resulting partition is $m$-regular.  Reverse by reading the $m$-ary expansion of the number of appearances of $j$, and turning $a_{j,k} \cdot m^k$ appearances of $j$ into $a_{j,k}$ appearances of $j m^k$.

\phantom{.}

\noindent \textbf{Example:} Consider the partition $\lambda = (108,18,18,18,18) \vdash 180$ as a 6-distinct partition.  Since $108 = 3 \cdot 6^2$, and $18 = 3 * 6^1$ appearing 4 times, we have $\phi_6(\lambda) = (3,\dots,3)$, where 3 appears $36 + 4 \cdot 6 = 60$ times.

\phantom{.}

Although Glaisher's map appears to be directional and applied to a proper subset of partitions, it can in fact be extended to an involution on all partitions in a natural way: construct a list of matrices $M_j$ indexed by the numbers $j$ not divisible by $m$.  Enumerate columns and rows starting with 0.  If the part $N = j m^k$ appears $a_{N,0} m^0 + a_{N,1} m^1 + a_{N,2} m^2 + \dots$ times as written in the ascending $m$-ary notation, then assign row $k$ of matrix $M_j$ to be $$\begin{matrix} a_{k,0} & a_{k,1} & a_{k,2} & \dots \end{matrix} \quad := \quad \begin{matrix} a_{N,0} & a_{N,1} & a_{N,2} & \dots \end{matrix}.$$   For example, if $m=2$, the partition $(20,5,5,4,2,2,1,1,1,1,1)$ of 43 would be depicted (with all other entries 0):

\begin{center} \begin{tabular}{ccc} $\begin{array}{c|ccc} 1 & & & \\ \hline & 1 & 0 & 1 \\ & 0 & 1 & 0 \\ & 1 & 0 & 0 \end{array}$ & $\begin{array}{c|ccc} 3 & & & \\ \hline & 0 & 0 & 0 \\ & 0 & 0 & 0 \\ & 0 & 0 & 0 \end{array}$ & $\begin{array}{c|ccc} 5 & & & \\ \hline & 0 & 1 & 0 \\ & 0 & 0 & 0 \\ & 1 & 0 & 0 \end{array}$ . \end{tabular} \end{center}

Now Glaisher's original map from $m$-distinct to $m$-regular partitions is simply transposition of these matrices applied to those in which only the first columns have nonzero entries.  Ignoring this restriction, we obtain a map $\phi_m$ which restricts to Glaisher's map on the relevant set, and which we shall therefore refer to freely as Glaisher's map.  Any other transformation of the matrices which preserves the sum in an antidiagonal is a weight-preserving map on partitions as well.  Of use in this paper will be the transformation on matrices with empty first column, which moves all other entries diagonally one space down and left and transposes the zero column to the top row.  (These \emph{part-frequency matrices} have seen use by several authors in recent years; see \cite{PFSurvey} for a recent survey.)

For the remainder of the paper, if we write partitions as parts with exponents, $(\lambda_1^{e_1}, \lambda_2^{e_2},\dots)$, we mean the partition in which part $\lambda_1$ appears $e_1$ times, part $\lambda_2$ appears $e_2$ times, and so forth.  The image partition in the map example would be written $(3^{60})$.

The goal of this paper is to consider the set of partitions which are simultaneously $s$-regular and $t$-distinct.  These were considered in an earlier paper of the author's \cite{Keith1}.  It is not difficult to see that the generating function for these partitions is $$\prod_{k=1}^\infty \frac{(1-q^{sk})(1-q^{tk})}{(1-q^k)(1-q^{stk})}.$$

From the fact that this generating function is symmetric in $s$ and $t$ we get immediately the following observation:

\begin{theorem} The number of partitions of $n$ which are simultaneously $s$-regular and $t$-distinct equals the number of partitions of $n$ which are simultaneously $t$-regular and $s$-distinct.
\end{theorem}

This observation cries out for a similarly explanatory bijection, which it is the purpose of this note to provide.

\phantom{.}

\noindent \textbf{Remark:} After first preprint distribution of this paper, Prof. Matja\v{z} Konvalinka kindly directed the author to a 1988 paper of Kathleen O'Hara, \cite{KOH}, which produces such a map.  The bijection of that paper applies to any two equivalent sets which are what Andrews calls \emph{partition ideals of order 1}, i.e. sets of partitions completely describable by listing maximum permissible numbers of repetitions of each part size independent of any others.  Certainly $s$-regular, $t$-distinct partitions (less than $t$ of any part size, less than 1 if divisible by $s$) are such a set.  Readers may be interested in that paper as well as in Konvalinka and Pak's analysis of its efficiency and geometry \cite{KP}.

The map in this paper is specifically designed for the sets under study and might be considered somewhat ``tuned'' compared to O'Hara's map. It is certainly distinct, since it will be shown later that the bijection is in fact an ensemble of possible bijections, depending on the order chosen for primes dividing $s$ and $t$.  As is common with combinatorial maps, it therefore gives rise to bijective proofs of numerous different intermediate results, some of which we list in the later part of the note.

\section{The bijection}

In \cite{Keith1}, it was established that a double use of Glaisher's map suffices in the case when $s$ and $t$ are coprime.

\begin{theorem} If $s$ and $t$ are coprime, then ${\phi_s} \phi_t$ maps $s$-regular, $t$-distinct partitions to $t$-regular, $s$-distinct partitions.
\end{theorem}

Here it is important that $\phi_t$ is performed first.  The result of this stage of the map is an $s$-regular, $t$-regular partition.  We observe that, \emph{in the case when $s$ and $t$ are coprime}, this set of partitions has the same generating function as our starting and ending sets.  However, this is \emph{not} the case when $s$ and $t$ are not coprime, for the set of $s$-regular, $t$-regular partitions has generating function $$\prod_{k=1}^\infty \frac{(1-q^{sk})(1-q^{tk})}{(1-q^k)(1-q^{\text{lcm}(\{s,t\})k})}$$

\noindent and $\text{lcm}(s,t) \neq st$ when $s$ and $t$ are not coprime.  Likewise, the set of $s$-distinct, $t$-distinct partitions is just the set of $\text{min}(s,t)$-distinct partitions.

A conjecture advanced in that paper, by a collaborator of the co-author, was that iteration of this map suffices to produce a bijection. Formally, 

\begin{conjecture}\label{DoubleGlaisher} There exists $\ell$, varying with $\lambda$, such that $({\phi_s} \phi_t)^\ell$ maps an $s$-regular, $t$-distinct partition $\lambda$ to a unique $s$-distinct, $t$-regular partition, with no intervening $(\phi_s \phi_t)^k$ being $s$-regular and $t$-distinct.
\end{conjecture}

Unfortunately, this conjecture fails, although there seems to be considerable interest remaining in exploring its domain of applicability and other properties.  We discuss these questions in the next section.

However, it is possible to give a bijection between the two sets.  This was achieved for 2-regular, $t$-distinct partitions in \cite{KeithMunagi}, where it was a useful tool in handling certain overpartitions, and the full result can now be announced.  The double-Glaisher map described above remains an important component.

\subsection{Preliminaries}

We desire to construct a map from the set of $s$-regular, $t$-distinct partitions to the set of $t$-regular, $s$-distinct partitions.  Let $\lambda$ be the partition to be mapped.

Let the prime factorizations of $s$ and $t$ be given by $$s = {p_1}^{e_1} {p_2}^{e_2} \dots {p_t}^{e_t} {q_1}^{d_1} \dots {q_r}^{d_r} \quad , \quad t = {p_1}^{b_1} {p_2}^{b_2} \dots {p_t}^{b_t} {v_1}^{c_1} \dots {v_u}^{c_u}.$$  All exponents are strictly positive, so the primes $p_i$ are common to both $s$ and $t$, and we name primes so that the primes $q_i$ and $v_i$ are exclusive to each modulus separately.

Collect the subset of parts of $\lambda$ which are not divisible by ${q_1}^{d_1} \dots {q_r}^{d_r}$.  This subset forms a $s^\prime := {q_1}^{d_1} \dots {q_r}^{d_r}$-regular, $t$-distinct partition, and $t$ is coprime with $s^\prime$.  Hence the double-Glaisher map $\phi_{s^\prime} \phi_t$ applied to this subpartition produces a $t$-regular, $s^\prime$-distinct partition.

We now perform an important ``wrapping'' step that will be repeated later.

The remaining parts are all divisible by $s^{\prime}$.  We divide all remaining parts of the partition by $s^\prime$.  After applying the remaining steps discussed below, to produce various subpartitions, we will multiply the resulting frequencies by $s^{\prime}$.  The resulting subpartition in the image can be separated from the image of the parts in the previous step by extracting the number of parts of each size modulo $s^{\prime}$.

We will proceed through the various primes $p_i$, ``wrapping'' successively with the powers ${p_i}^{e_i}$ after each step.  In each case, the subpartitions so constructed can be separated from the rest of the image partition by noting whether their frequency of appearance is divisible by ${p_i}^{e_i}$.

The choice of ordering on the primes in $s$ does matter, as we will show by example later, meaning many possible maps can be constructed.  Any, however, suffice as a bijectiion.

\subsection{The map for each prime $p_i$}

Begin with $p_1$.  Write $t = {p_1}^{b_1} k$.  Collect the set of parts remaining in $\lambda$ that are not divisible by ${p_1}^{e_1}$.  These form a ${p_1}^{e_1}$-regular, $t$-distinct partition.  Call it $\rho$.

Suppose that a part size $m$ appears $a + Ck$ times in $\rho$, $0 \leq a < k$, $0 \leq C < {p_1}^{b_1}$.

Collect the $a$ appearances of $m$ over all part sizes.  The collection of these constitutes a partition which is ${p_1}^{e_1}$-regular and $k$-distinct, with $k$ coprime to ${p_1}^{e_1}$.  The bijection $\phi_{{p_1}^{e_1}} \phi_k$ can thus be applied.  The result is a set of parts forming a partition which is $k$-regular and ${p_1}^{e_1}$-distinct.

Now consider separately the parts appearing $Ck$ times.  Divide the frequency of appearance of each part size $x$ by $k$ and multiply the size of each part by $k$.  The result is a collection of parts $xk$ divisible by $k$ but not by ${p_1}^{e_1}$ (since $k$ and ${p_1}^{e_1}$ are coprime), and appearing fewer than ${p_1}^{b_1}$ times.  Apply $\phi_{p_1}$ to obtain a collection of parts divisible by $k$ but not by $k {p_1}^{b_1} = t$, and appearing fewer than ${p_1}^{e_1}$ times.

Together with the previous collection, we have a $t$-regular, ${p_1}^{e_1}$-distinct partition.  We conclude the ``wrapping'' step by multiplying frequencies by any $s^\prime$.

\phantom{.}

\noindent \textbf{Example:} Consider $\lambda = (10,10,10,10,5,5,5,5,5,5,5,3,3,3,3,3,1,1)$ as a 9-regular, 15-distinct partition.  We have $t = 3^1 5$, $r=1$, $k=5$, and no $s^\prime$ wrapping step.

\begin{itemize}
\item The part size 10 appears 4 times: $a=4$, $C=0$.  Set aside the subpartition $(10,10,10,10)$.
\item The part size 5 appears 7 times: $a=2$, $C=1$.
\begin{itemize}
\item Set aside the subpartition $(5,5)$.
\item For the subpartition $(5,5,5,5,5)$, we have $x=5$, $Ck = 5 \cdot 1$.  Thus this is replaced by 1 appearance of $5 \cdot 5 = 25$.  The partition $(25)$ is fixed under $\phi_3$.
\end{itemize}
\item The part size 3 appears 5 times: $a=0$, $C=1$.
\begin{itemize}
\item For the subpartition $(3,3,3,3,3)$, we replace it by part 15, then, applying $\phi_3$, obtain $(5,5,5)$.
\end{itemize}
\item The part size 1 appears 2 times.  Set aside the subpartition $(1,1)$.
\item Taking the partition $(10,10,10,10,5,5,1,1)$ as a 3-regular, 5-distinct partition, we employ double-Glaisher ${\phi_3} \phi_5$ and obtain the 5-regular, 3-distinct partition $(18,18, 9, 3, 2, 2)$.
\end{itemize}

Taking the union of all parts, the image partition is $(25,18,18,9,5,5,5,3,2,2)$, a 15-regular, 9-distinct partition.

\phantom{.}

\subsection{Iteration over primes}

The rest of the map is simple iteration. The remaining parts are all divisible by ${p_1}^{e_1}$; wrap with ${p_1}^{e_1}$ and repeat.  That is, divide all remaining parts by ${p_1}^{e_1}$, pick another prime, perform the map in the previous subsection, and multiply the frequencies of the resulting partition by ${p_1}^{e_1}$, followed by $s^\prime$ if necessary.  Continue nesting any wrappings required.

\phantom{.}

\noindent \textbf{Remark:} To see that the order in which primes are chosen matters, consider the partition $(9)$ as an 18-regular, 30-distinct partition.  Since $s^\prime = 1$ there is no initial wrapping step.  If $p_1 = 3$, then the part of size 9 is handled as $(1)$ under wrapping by 9, and its image is $(1^9)$; if $p_1 = 2$, then $(9)$ is a 2-regular, 15-distinct partition and is fixed by the map.

\phantom{.}

\noindent \textbf{Proof of injectivity:} The parts produced at every step are $t$-regular, and the wrapping steps multiply frequencies, so it remains to show that the resulting partition is $s$-distinct.  Bearing in mind the remark above, assume that we have a preferred order on the primes $p_i$.

Frequencies of parts contributed in the first step are less than $s^\prime$.  Frequencies of parts contributed in the next step are divisible by $s^\prime$, and are strictly less than ${p_1}^{e_1} s^\prime$.  The sum of the two frequencies, should a part size have appeared at both steps, is any possible value from 0 to ${p_1}^{e_1} s^\prime - 1$.  Frequencies of parts contributed in the next step are divisible by ${p_1}^{e_1} s^\prime$, and are strictly less than ${p_1}^{e_1} {p_2}^{e_2} s^\prime$.  The sum of these with the previous frequencies is any possible value up to ${p_1}^{e_1} {p_2}^{e_2} s^\prime - 1$.  This continues up to $s$ in the last step.

\phantom{.}

\noindent \textbf{Proof of reversibility:} To establish that the map is a bijection it suffices to confirm reversibility.

All steps after the initial $s^\prime$ step had their frequencies multiplied by $s^\prime$.  Hence, the parts which appeared at the first step were necessarily the parts contributing to any nonzero frequency mod $s^\prime$.  These can be extracted and, forming a $t$-regular, $s^\prime$-distinct partition, the double-Glaisher map reverses the initial step, yielding the subpartition which is $t$-distinct and $s^\prime$-regular.

Of the remaining parts, all contributions after the next step had their frequencies multiplied by ${p_1}^{e_1} s^\prime$.  We extract the nonzero portions of any remaining frequencies of appearance mod ${p_1}^{e_1} s^\prime$  Dividing the frequency of appearance of each part size by $s^\prime$, we obtain a $t$-regular, ${p_1}^{e_1}$-distinct subpartition.   We apply $\phi_{p_1}$ and separately consider resulting parts divisible by $k$ (recall that $t = {p_1}^{b_1} k$) and not; the subpartition consisting of those not divisible by $k$ have $\phi_k$ applied, and those divisible by $k$ have their sizes divided by $k$ and their frequencies multiplied by $k$.  Finally, we finish ``uunwrapping'' by multiplying all resulting sizes by $s^\prime$.  This reverses the map at this step.  The resulting subpartition is $t$-distinct and ${p_1}^{e_1} s^\prime$-regular, with all part sizes divisible by $s^\prime$ and thus distinguishable from those arising in the previous step.

Repeat for all remaining primes in $s$, and reversibility is shown. \hfill $\Box$

\phantom{.}

\subsection{A visualization.} Glaisher's original bijection can be viewed as a matrix transposition; O'Hara's algorithm has a geometric interpretation studied by Konvalinka and Pak in \cite{KP}.  The algorithm of this paper, as well, can be visualized in a fashion somewhat intermediate to the two, and which may be of interest to readers.

The first step of the algorithm, the $\phi_t \phi_{s^\prime}$ step, can be thought of as two of Glaisher's matrix transpositions, with the necessity of rewriting the matrices in the two different bases between steps.  The ${p_i}^{e_i}$ steps can be visualized similarly, with the $a$ and the $Ck$ slices being treated separately.

Suppose that we have performed the division by $s^\prime$, and are now handling the parts of the wrapped partition not divisible by ${p_1}^{e_1}$.  We now extend the matrix concept to a three-dimensional array: consider arrays labeled by $m \not\equiv 0 \pmod{p_1}$, and $e_1$ columns labeled $m$ through $m{p_1}^{e_1-1}$, with each row in the $m{p_1}^i$ column being associated to part size $m {p_1}^i k^j$ for some $j$.  Each column extends indefinitely far in the $k^j$ direction.  In row $m {p_1}^i k^j$ write out the $k$-ary notation for number of appearance of the part size.  The $0 \leq a < k$ part constitutes one face of the prism; the $Ck$ part the remainder, with $C < {p_1}^{b_1}$.  This is illustrated below.

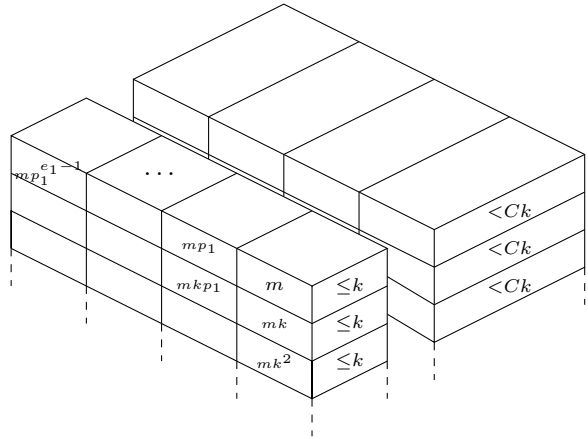
\begin{center}\begin{tabular}{c}
\begin{tikzpicture}[scale=0.5]
\draw (0,-1) -- (0,-4) -- (0,-3) -- (-8,1) -- (-8,0) -- (-8,3) -- (0,-1) -- (2,0) -- (-6,4) -- (-8,3);
\draw[dashed] (0,-4) -- (0,-5);
\draw[dashed] (2,-3) -- (2,-4);
\draw[dashed] (-2,-3) -- (-2,-4);
\draw[dashed] (-4,-2) -- (-4,-3);
\draw[dashed] (-6,-1) -- (-6,-2);
\draw[dashed] (-8,0) -- (-8,-1);
\draw (-4,3) -- (-6,2) -- (-6,-1);
\draw (-2,2) -- (-4,1) -- (-4,-2);
\draw (0,1) -- (-2,0) -- (-2,-3);
\draw (-8,2) -- (0,-2) -- (2,-1);
\draw (0,-3) -- (2,-2) -- (2,0);
\draw (2,-2) -- (2,-3) -- (0,-4) -- (-8,0);
\draw (-1,-1) node{$\scriptstyle m$};
\draw (1,-1) node{$\scriptstyle \leq k$};
\draw (-1,-2) node{$\scriptscriptstyle mk$};
\draw (1,-2) node{$\scriptstyle \leq k$};
\draw (-1,-3) node{$\scriptscriptstyle mk^2$};
\draw (1,-3) node{$\scriptstyle \leq k$};
\draw (-3,0) node{$\scriptscriptstyle mp_1$};
\draw (-3,-1) node{$\scriptscriptstyle mkp_1$};
\draw (-4,2) node{$\dots$};
\draw (-7,2) node{${\scriptscriptstyle mp_1^{e_1-1}}$};
\draw (-4.75,3.375) -- (-4.75,4.5) -- (3.25,0.5) -- (3.25,-2.5) -- (2,-1.75);
\draw (3.25,-1.5) -- (2,-0.75);
\draw (3.25,-0.5) -- (-4.75,3.5);
\draw (-2.75,2.375) -- (-2.75,3.5);
\draw (-0.75,1.375) -- (-0.75,2.5);
\draw (1.25,0.375) -- (1.25,1.5);
\draw (-4.75,4.5) -- (-0.75,6.5) -- (7.25,2.5) -- (7.25,-0.5) -- (3.25,-2.5);
\draw (3.25,-1.5) -- (7.25,0.5);
\draw (3.25,-0.5) -- (7.25,1.5);
\draw (3.25,0.5) -- (7.25,2.5);
\draw (1.25,1.5) -- (5.25,3.5);
\draw (-0.75,2.5) -- (3.25,4.5);
\draw (-2.75,3.5) -- (1.25,5.5);
\draw[dashed] (7.25,-0.5) -- (7.25,-1.5);
\draw (5.25,-1) node{$\scriptstyle <Ck$};
\draw (5.25,0) node{$\scriptstyle <Ck$};
\draw (5.25,1) node{$\scriptstyle <Ck$};
\draw[dashed] (3.25,-2.5) -- (3.25,-3.5);
%\foreach \x in {0,1,2,3,4,5}
% {\draw (\x,0.5) node{\x};
%  \draw (-0.5,-\x) node{\x};
%  }
\end{tikzpicture} \\
The initial partition: size axes left and down, frequency axis in base $k$ to the right.
\end{tabular}
\end{center}

In the base $k$, our operations now are the following: the $a$ face is transposed around the $m {p_1}^i$ axis to become the top face of the prism, while the $Ck$ portion is shifted diagonally down and back one space.

The resulting occupied shape is then rewritten in base $p_1$ and the entire shape transposed around the $m k^j$ axis.  The result is parts which, if divisible by $k$, are not divisible by ${p_1}^{b_1}$, and if divisible by ${p_1}^{b_1}$, are not divisible by $k$; all appear less than ${p_1}^{e_1}$ times.  The resulting visualization is illustrated below.

\begin{center}\begin{tabular}{cc}
\begin{tikzpicture}[scale=0.5]
\draw (-2,4) -- (-4,3) -- (-4,2) -- (0,0) -- (2,1) -- (2,-2);
\draw (-4,3) -- (0,1) -- (2,2);
\draw (0,1) -- (0,-3);
\draw (-4,2) -- (-4,-1);
\draw[dashed] (-4,-1) -- (-4,-2);
\draw[dashed] (0,-3) -- (0,-4);
\draw[dashed] (2,-2) -- (2,-3);
\draw[dashed] (2,1) -- (3,1.5);
\draw[dashed] (2,2) -- (3,2.5);
\draw[dashed] (-2,4) -- (-1,4.5);
\draw (-0.5,0.75) node{$\scriptstyle m$};
\draw (-1,-0.5) -- (-1,1.5);
\draw (-0.5,-0.25) node{$\scriptscriptstyle mk$};
\draw (-0.5,-1.25) node{$\scriptscriptstyle \dots$};
\draw (-1.5,1.25) node{$\scriptstyle \dots$};
\draw (0.75,-0.25) node{$\scriptscriptstyle <{p_1}^{b_1}$};
\draw (2,0) -- (0,-1) -- (-4,1);
\draw (-3,2.5) -- (-3,0.5);
\draw (-5.5,2.25) node{$\scriptscriptstyle m{p_1}^{e_1-1} \rightarrow$};
\draw (-5.5,1.25) node{$\scriptscriptstyle mk{p_1}^{e_1-1} \rightarrow$};
\end{tikzpicture}
&
\begin{tikzpicture}[scale=0.5]
\draw (2,4) -- (4,3) -- (4,2) -- (0,0) -- (-2,1) -- (-2,-2);
\draw (4,3) -- (0,1) -- (-2,2);
\draw (0,1) -- (0,-3);
\draw (4,2) -- (4,-1);
\draw[dashed] (4,-1) -- (4,-2);
\draw[dashed] (0,-3) -- (0,-4);
\draw[dashed] (-2,-2) -- (-2,-3);
\draw[dashed] (-2,1) -- (-3,1.5);
\draw[dashed] (-2,2) -- (-3,2.5);
\draw[dashed] (2,4) -- (1,4.5);
\draw (1,-0.5) -- (1,1.5);
\draw (-2,0) -- (0,-1) -- (4,1);
\draw (3,2.5) -- (3,0.5);
\end{tikzpicture}
\\
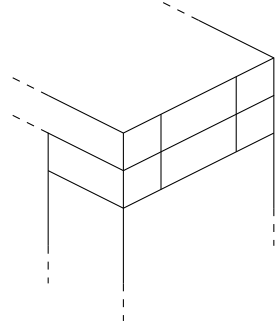
After first transposition. & Rewrite frequency axis in base $p_1$, and transpose left-right.
\end{tabular}
\end{center}

\section{Additional questions; the double-Glaisher conjecture}

In combinatorics it is often the case that a bijection permits the exploration of a multitude of statistics that are naturally preserved by the map.  It is to be hoped that this map will yield such results as well. For instance, one may observe the following.  Defining $s^\prime$ as in the previous section, the following theorem is easy.

\begin{theorem} The set of $s$-regular, $t$-distinct partitions of $n$ in which not all parts are divisible by $s^\prime$ is equinumerous with the set of $t$-regular, $s$-distinct partitions of $n$ in which not all frequencies are divisible by $s^\prime$.
\end{theorem}

We also obtain any intermediate sets as equivalent sets.  For instance, after the first step, we obtain the following.  (Here scalar multiplication of a partition $\lambda \vdash n$ by $c$ means to multiply all sizes by $c$, obtaining a partition $c \lambda \vdash cn$.)

\begin{theorem} The set of $s$-regular, $t$-distinct partitions of $n$ is equinumerous with the set of pairs of partitions $(\lambda,\mu)$ in which $\mu$ is $t$-regular and $s^\prime$-distinct, and $\lambda = s^\prime \gamma$, with $\gamma$ being $s/s^\prime$-regular and $t$-distinct, in which $\vert \lambda \vert + \vert \mu \vert = n$.
\end{theorem}

After the $p_1$ step we have the following equivalence.

\begin{theorem} The set of $s$-regular, $t$-distinct partitions of $n$ is equinumerous with the set of pairs of partitions $(\lambda,\mu)$ in which $\mu$ is $t$-regular and ${p_1}^{e_1} s^\prime$-distinct, and $\lambda = {p_1}^{e_1} s^\prime \gamma$, with $\gamma$ being $s/(s^\prime {p_1}^{e_1})$-regular and $t$-distinct, in which $\vert \lambda \vert + \vert \mu \vert = n$.
\end{theorem}

Of course, if one took only a subset of the factors of $s^\prime$ and went partway through each $p_i$ step, one could obtain an intermediate theorem of this type for \emph{any} divisor $d \vert s$.  To shorten notation let $f_k = \prod_{i=1}^\infty (1-q^{ik})$.  Each of these intermediate theorems is simply a map showing the equivalence of the generating functions $$\frac{f_t f_d}{f_1 f_{td}} \cdot \frac{f_{dt} f_s}{f_d f_{st}} = \frac{f_t f_s}{f_1 f_{st}}.$$

\phantom{.}

\subsection{A counterexample to Conjecture \ref{DoubleGlaisher}.} An example where the map of Conjecture \ref{DoubleGlaisher} works would be the partition $\lambda = (50)$ considered as a 6-regular, 10-distinct partition.  We obtain the following chain of images:

\begin{align*}
\phi_{10}(\lambda) &= (5,5,5,5,5,5,5,5,5,5) \, ;& \, {\phi_6}\phi_{10}(\lambda) &= (30,5,5,5,5) =: \beta \, ; \\
\phi_{10}(\beta) &= (5,5,5,5,3,\dots,3) \, ;& \, {\phi_6}\phi_{10}(\beta) &= (18,5,5,5,5,3,3,3,3) .
\end{align*}

The final partition is the first in the series of images to be 10-regular and 6-distinct.  (Although note that the conjecture would not be disproved if an image halfway through the map, after $\phi_{10}$, happened to satisfy these conditions.)

However, a counterexample can be found.  Consider $\lambda = (108,18,18,18,18)$ as a 10-regular, 6-distinct partition.  We have $${\phi_{10}} \phi_{6} (\lambda) = (30,30,30,30,30,30).$$  This is not 6-distinct and so we cannot apply $\phi_6$ as originally defined, and it is a fixed point of $\phi_6$ considered as matrix transposition.  The sequence now starts working backward, and $(108,18^4)$ is a fixed point of $\phi_{10}$.

The sequence of partitions obtained by the generalized map is thus as follows.

\begin{align*}
(108,18^4) & {{\phi_6} \atop \longrightarrow} & (3^{60}) \quad \quad & {{\phi_{10}} \atop \longrightarrow} (30^6) \\
& {{\phi_6} \atop \longrightarrow} & (30^6) \quad \quad & {{\phi_{10}} \atop \longrightarrow} (3^{60}) \\
& {{\phi_6} \atop \longrightarrow} & (108,18^4) \quad \quad & {{\phi_{10}} \atop \longrightarrow} (108,18^4)
\end{align*}

Thus $\phi_{10} \phi_6$ does not map $(108,18^4)$, a 10-regular, 6-distinct partition, to a 6-regular, 10-distinct partition after any number of iterations: $(30^6)$ is not 6-regular, $(3^{60})$ is not 10-distinct, and $(108,18^4)$ is not 6-regular.

Nevertheless, Conjecture \ref{DoubleGlaisher} remains interesting for several reasons.  Most importantly, it clearly works for some partitions.  Therefore, some natural questions are:

\phantom{.}

\noindent \textbf{Question:} What is the domain of applicability of the recipe $(\phi_s \phi_t)^\ell$ for mapping $s$-regular, $t$-distinct to $t$-regular, $s$-distinct partitions - i.e., for what sets of partitions can we guarantee that the method gives a bijection?

\phantom{.}

\noindent \textbf{Question:} For such a subset of partitions, what is the distribution of the necessary $\ell$?

\phantom{.}

\noindent \textbf{Question:} What is the distribution of the lengths of the cycles when they occur?

\phantom{.}

For example, the partitions which take 0 iterations are just exactly the set of partitions which are in both sets: if $s < t$, then these are the partitions which are simultaneously $s$-regular, $t$-regular, and $s$-distinct.  These partitions have received some recent study when $s=2$ in \cite{DremaSaikia}, albeit from a different point of view.  These partitions are counted by the following generating function.

\begin{theorem} Say $s < t$.  The generating function for the number of partitions of $n$ which are simultaneously $s$-regular, $t$-regular and $s$-distinct is $$\frac{f_s f_t}{f_1 f_{\text{lcm}(s,t)}} \cdot \frac{f_s f_{s \cdot \text{lcm}(s,t)}}{f_{s^2} f_{st}}.$$
\end{theorem}

What then are the generating functions for the subset that requires 1 iteration, 2 iterations, etc.?

\phantom{.}

Glaisher's map is central to many combinatorial arguments in partitions.  Hopefully, this map, while necessarily somewhat more complicated, will also be of use to investigators of the subject.

\phantom{.}

\noindent \textbf{Acknowledgement:} This map was originally presented in the Specialty Seminar in Partition Theory, $q$-Series and Related Topics hosted by Michigan Technological University in the Spring of 2022.  Several steps were combined and overall explication improved in the course of work with undergraduate student Hadley Wells as part of a research project in the summer of 2022.  The author also thanks Prof. Konvalinka for the communication mentioned earlier.

\end{document}